\documentclass{amsart}
\usepackage{amsmath}
\usepackage{amssymb}
\usepackage{amsthm}
\usepackage{epsfig}
\usepackage{amsfonts}
\usepackage{amscd}
\usepackage{mathrsfs}
\usepackage{enumerate}
\usepackage{latexsym}
\usepackage{url}

\newtheorem{theorem}{Theorem}[section]
\newtheorem{lemma}[theorem]{Lemma}

\newtheorem{notationandlemma}[theorem]{Notation and Lemma}
\newtheorem*{xclaim}{Claim}
\newtheorem*{xconjecture}{Conjecture}

\newenvironment{proofoftheclaim}{\noindent\textit{Proof of the claim.}}{\hfill{}\newline}
\newenvironment{proof*}{\noindent\textbf{Proof.}}{\hfill{\;$\blacksquare$}\newline}

\newenvironment{proof of the claim}{\noindent\textbf{Proof of the claim.}}{\hfill{$\square$}\newline}

\theoremstyle{definition}

\newtheorem{example}[theorem]{Example}

\newtheorem*{xremark}{Remark}

\theoremstyle{remark}

\numberwithin{equation}{section}

\begin{document}

\title[One-point connectifications]{The existence of one-point connectifications}

\author{M.R. Koushesh}

\address{Department of Mathematical Sciences, Isfahan University of Technology, Isfahan 84156--83111, Iran, and,}

\address{School of Mathematics, Institute for Research in Fundamental Sciences (IPM), P.O. Box: 19395--5746, Tehran, Iran.}

\email{koushesh@cc.iut.ac.ir}

\thanks{This research was in part supported by a grant from IPM (Nos. 95030418 and 96030416).}




\begin{abstract}
  P. Alexandroff proved that a locally compact $T_2$-space has a $T_2$ one-point compactification (obtained by adding a ``point at infinity'') if and only if it is non-compact. He also asked for characterizations of spaces which have one-point connectifications. Here, we study one-point connectifications, and in analogy with Alexandroff's theorem, we prove that in the realm of $T_i$-spaces ($i=3\frac{1}{2},4,5$) a locally connected space has a one-point connectification if and only if it has no compact component. We extend this theorem to the case $i=6$ by assuming some set-theoretic assumption, and to the case $i=2$ by slightly modifying its statement. We further extended the theorem by proving that a locally connected metrizable (resp. paracompact) space has a metrizable (resp. paracompact) one-point connectification if and only if it has no compact component. Contrary to the case of the one-point compactification, a one-point connectification, if exists, may not be unique. We consider the collection of all one-point connectifications of a locally connected locally compact space in the realm of $T_i$-spaces ($i=3\frac{1}{2},4,5$). We prove that this collection, naturally partially ordered, is a compact conditionally complete lattice whose order-structure determines the topology of all Stone--\v{C}ech remainders of components of the space.
\end{abstract}

\maketitle

\tableofcontents

\section{Introduction}

A \textit{one-point connectification} (\textit{one-point compactification}, respectively) of a space $X$ is a connected (compact, respectively) space $Y$ which contains $X$ as a dense subspace and $Y\setminus X$ is a singleton.

In 1924, P. Alexandroff proved in \cite{A} that a locally compact non-compact Hausdorff space $X$ can be compactified by adding a ``point at infinity.'' The analogous question (also attributed to Alexandroff) of determining spaces which can be connectified by adding one point remained however unsettled. The question indeed has motivated a significant amount of research, with perhaps the earliest serious work dating as far back as 1945. In 1945, B. Knaster proved in \cite{Kn} that a separable metrizable space has a separable metrizable one-point connectification if and only if it can be embedded in a connected separable metrizable space as a proper open subspace. Among more recent results we mention that of M. Abry, J.J. Dijkstra, and J. van Mill in \cite{ADvM} stating that a separable metrizable space in which every component is open has a separable metrizable one-point connectification if and only if it has no compact component. Other results concerning connectifications of spaces may be found in \cite{ATTW}, \cite{DvM}, \cite{DvMV}, \cite{DV}, \cite{DT}, \cite{GKLD}, \cite{HB}, \cite{Ko}, \cite{PW1} and \cite{WW}. For results concerning general one-point extensions of spaces see \cite{Ko5} and \cite{MRW0} and the cited bibliography therein.

Here, we start with Alexandroff's one-point compactification theorem, which suitably reformulated, states that in the realm of $T_i$-spaces ($i=2,3,3\frac{1}{2},4,5$) a locally compact space has a one-point compactification if and only if it is non-compact. Motivated by this, we prove that in the realm of $T_i$-spaces ($i=3\frac{1}{2},4,5$) a locally connected space has a one-point connectification if and only if it has no compact component. We extend this theorem to the case $i=6$ by assuming certain set-theoretic assumption, and to the case $i=2$ by slightly modifying its statement, proving that a locally connected $T_2$-space has a $T_2$ one-point connectification if and only if it has no almost compact component. We further extend the theorem by proving it in the realm of $T_i$-spaces ($i=3\frac{1}{2},4,5$) having a ceratin topological property $\mathscr{P}$. The topological property $\mathscr{P}$ is subject to some mild requirements and includes most ``covering properties'' (i.e., topological properties defined by means of existence of certain kinds of open subcovers or refinements of open covers of certain types) among them are realcompactness, pseudocompactness, paracompactness, metacompactness, and the Lindel\"{o}f property. (Indeed, in the case when $\mathscr{P}$ is either  realcompactness or pseudocompactness our theorem holds even in the realm of $T_6$-spaces, assuming some set-theoretic assumption.) In particular, we prove that a locally connected metrizable space has a metrizable one-point connectification if and only if it has no compact component. Contrary to the case of the one-point compactification, a one-point connectification, if exists, may not be unique. We consider the collection of all one-point connectifications of a locally connected locally compact space in the realm of $T_i$-spaces ($i=3\frac{1}{2},4,5$). We prove that this collection, naturally partially ordered, is a compact conditionally complete lattice whose order-structure determines the topology of all Stone--\v{C}ech remainders of components of the space.

In what follows we will use ceratin basic facts from the theory of the Stone--\v{C}ech compactification. We review some of these facts in the following and refer the reader to the texts \cite{GJ} and \cite{PW} for further information.

The \textit{Stone--\v{C}ech compactification} of a completely regular space $X$, denoted by $\beta X$, is the (unique) compactification of $X$ with the property that every continuous mapping $f:X\rightarrow[0,1]$ extends to a continuous mapping $F:X\rightarrow[0,1]$. For a completely regular space $X$ the Stone--\v{C}ech compactification of $X$ always exists.

For a completely regular space $X$ the Stone--\v{C}ech compactification has the following properties:
\begin{itemize}
  \item $X$ is locally compact if and only if $X$ is open in $\beta X$.
  \item $\beta U=\mathrm{cl}_{\beta X}U$ for every open and closed subspace $U$ of $X$.
  \item $\beta T=\beta X$ for every subspace $T$ of $\beta X$ containing $X$.
  \item Disjoint open and closed subspaces of $X$ have disjoint closures in $\beta X$.
\end{itemize}

\textit{Components} (also called \textit{connected components}) of a space $X$ are maximal connected subspaces of $X$ which partition $X$ into pairwise disjoint non-empty closed subspaces. Components of locally connected spaces are open. Thus, in particular, when $X$ is a locally connected completely regular space, components are open and closed subspaces in $X$, and therefore, their closures in $\beta X$ form a collection of pairwise disjoint open and closed subspaces of $\beta X$.

A \textit{zero-set} in a space $X$ is a set of the form $f^{-1}(0)$, denoted by $Z(f)$, where $f:X\rightarrow[0,1]$ is a continuous mapping. It is known that a completely regular space $X$ is pseudocompct if and only if $\beta X$ has no non-empty zero-set disjoint from $X$. (A space $X$ is \textit{pseudocompct} if there is no unbounded continuous mapping $f:X\rightarrow\mathbb{R}$.) Observe that pseudocompactness coincides with countable compactness in the class of normal spaces. Therefore, a normal space $X$ is countably compact if and only if $\beta X$ has no non-empty zero-set disjoint from $X$.

\section{The existence of one-point connectifications of a space}

In this section we prove a theorem which provides necessary and sufficient conditions for a locally connected space to have a one-point connectification. This will be done in the classes of $T_i$-spaces where $i=3\frac{1}{2},4,5$. We extend this theorem to the case $i=6$ by assuming some set-theoretic assumption, and to the case $i=2$ by slightly modifying its statement. We further extend the theorem to other classes of spaces such as paracompact, metrizable, realcompact and pseudocompact spaces. The theorems may be all considered as being duals to the Alexandroff one-point compactification theorem.

In the following, we first review the statement of the Alexandroff one-point compactification theorem and we then see how the statement of this theorem (suitably reformulated for our purpose) motivates the formulation of the statements of our theorems.

\begin{theorem}[Alexandroff]\label{KLFA}
A locally compact Hausdorff space has a Hausdorff one-point compactification if and only if it is non-compact.
\end{theorem}

Indeed, for a locally compact non-compact Hausdorff ($T_2$) space $X$ a one-point compactification is a set $Y=X\cup\{p\}$, where $p$ is not in $X$, topologized by defining the open subspaces to be all the open subspaces of $X$ together with all sets of the form $Y\setminus C$ where $C$ is a compact subspace of $X$.

Observe that locally compact Hausdorff spaces are completely regular ($T_{3\frac{1}{2}}$) and compact Hausdorff spaces are normal ($T_4$). The Alexandroff theorem now motivates the following more general formulation.

Recall that a space is called \textit{hereditarily normal} (\textit{completely normal} or $T_5$) if all its subspaces are normal. It is known that a $T_1$-space is hereditarily normal if and only if every two separated subspaces are contained in disjoint open subspaces. (See Theorem 2.1.7 of \cite{E}. Two subspaces of a space are called \textit{separated} if neither intersects the closure of the other.)

\begin{theorem}\label{HFD}
Let $i=2,3,3\frac{1}{2},4,5$. A locally compact $T_i$-space has a $T_i$ one-point compactification if and only if it is non-compact.
\end{theorem}

\begin{proof*}
As pointed out earlier, locally compact Hausdorff spaces are completely regular and compact Hausdorff spaces are normal. We therefore need to check the theorem only in the case when $i=5$.

Let $X$ be a locally compact non-compact hereditarily normal space. Let $Y=X\cup\{p\}$ be the one-point compactification of $X$. Let $A$ and $B$ be a pair of separated subspaces in $Y$. We fined disjoint open subspaces in $Y$ containing $A$ and $B$.

Suppose that $p$ is in neither $A$ nor $B$. Then $A$ and $B$ are separated subspaces in $X$. Therefore, since $X$ is a hereditarily normal space, there are disjoint open subspaces in $X$ (which are also open in $Y$) containing $A$ and $B$.

Suppose that $p$ is in either $A$ or $B$, say in $A$. Then $p$ is not in $\mathrm{cl}_YB$ (as the latter is disjoint from $A$). Since $Y$ is regular, there is an open neighborhood $W$ of $p$ in $Y$ whose closure $\mathrm{cl}_YW$ is disjoint from $B$. Now $A\cap X$ and $B$ are separated subspaces in $X$. Since $X$ is a hereditarily normal space, there are disjoint open subspaces $U$ and $V$ in $X$ (which are also open in $Y$) containing $A\cap X$ and $B$, respectively. Now $U\cup W$ and $V\cap (Y\setminus\mathrm{cl}_YW)$ are disjoint open subspaces in $Y$ containing $A$ and $B$, respectively.
\end{proof*}

Observe that the above theorem does not hold in the case when $i=6$. (Recall that a space is called $T_6$, or \textit{perfectly normal}, if it is a normal space in which every closed subspace is a $G_\delta$-set.) To see this, let $X$ be an uncountable discrete space and let $Y=X\cup\{p\}$ be the one-point compactification of $X$. Any open subspace in $Y$ which has $p$ contain $X$ except for a finite number of its elements. Therefore, $\{p\}$ is a closed subspace in $Y$ which is not a $G_\delta$-set in $Y$.

Our next theorem is analogous to Theorem \ref{HFD}, replacing ``compactness'' by ``connectedness,'' with the exclusion of the cases when $i=2,3$ and the inclusion of the case when $i=6$. (The case when $i=2$ will be dealt with separately in Theorem \ref{JLJH}.) In the case when $i=6$ we assume Martin's Axiom ($\mathbf{MA}$) and the negation of the Continuum Hypothesis ($\mathbf{CH}$).

It is known that if the Continuum Hypothesis is negated and Martin's Axiom is assumed, then every perfectly normal, countably compact space is compact, or, equivalently, every perfectly normal non-compact space is not countably compact. (See \cite{W}.) This particularly implies that, under $\mathbf{MA}+\neg\mathbf{CH}$, any non-compact perfectly normal space $X$ is not countably compact, and therefore, there is a non-empty zero-set of $\beta X$ which is disjoint from $X$.

In the proof of the following theorem we will make use of some of the ideas we have already used in the proof of Theorem 2.2 in \cite{Ko}.

\begin{theorem}\label{PHRS}
Let $i=3\frac{1}{2},4,5$ or $i=6$ with $\mathbf{MA}+\neg\mathbf{CH}$. A locally connected $T_i$-space has a $T_i$ one-point connectification if and only if it has no compact component.
\end{theorem}

\begin{proof*}
Let $X$ be a locally connected $T_i$-space.

Suppose that $X$ has a $T_i$ one-point connectification $Y$. Suppose to the contrary that $X$ has a compact component $C$. Then $C$ is closed in $Y$, as $C$ is compact and $Y$ is Hausdorff. On the other hand, $C$ is also open in $Y$, as $C$ is open in $X$, since $X$ is locally connected, and $X$ is open in $Y$, since $Y$ is a $T_1$-space.
Therefore $C$ is a non-empty proper simultaneously open and closed subspace of $Y$, which contradicts connectedness of $Y$.

Suppose that $X$ has no compact component. We prove that $X$ has a $T_i$ one-point connectification. We divide the proof into verification of several claims.

\begin{xclaim}
Let $i=3\frac{1}{2},4,5$. There is a non-empty compact subspace $Z$ of $\beta X$ disjoint from $X$ which intersects the closure in $\beta X$ of every component of $X$.
\end{xclaim}

\begin{proofoftheclaim}
Let $C$ be a component of $X$. Then $C$ is non-compact and therefore $\mathrm{cl}_{\beta X}C\setminus X$ is non-empty. Let $z_C$ be in $\mathrm{cl}_{\beta X}C\setminus X$. Let
\[Z=\bigg(\bigcup_{\mbox{$C$ is a component of $X$}}\{z_C\}\bigg)\cup\bigg(\beta X\setminus\bigcup_{\mbox{$C$ is a component of $X$}}\mathrm{cl}_{\beta X}C\bigg).\]
Note that components of a space are always closed in the space. Furthermore, since $X$ is locally connected, components of $X$ are also open in $X$. Therefore, components of $X$, being open and closed subspaces of $X$, have open and closed closures in $\beta X$. Also, distinct components of $X$, being disjoint open and closed subspaces of $X$, have disjoint closures in $\beta X$. This, in particular, implies that $Z$ has no accumulation point in $\mathrm{cl}_{\beta X}C$ for every component $C$ of $X$, and is thus closed in $\beta X$ and is therefore compact. It is clear that $Z$ is non-empty, is disjoint from $X$, as is disjoint from every component of $X$, and intersects the closure in $\beta X$ of every component of $X$.

This completes the proof of the claim.
\end{proofoftheclaim}

\begin{xclaim}
Let $i=6$ $($and assume $\mathbf{MA}+\neg\mathbf{CH}$$)$. There is a non-empty compact $G_\delta$ subspace $Z$ of $\beta X$ disjoint from $X$ which intersects the closure in $\beta X$ of every component of $X$.
\end{xclaim}

\begin{proofoftheclaim}
Let $C$ be a component of $X$. Then $C$ is non-compact. Therefore, since $C$ is perfectly normal (as $X$ is perfectly normal and perfect normality is hereditary with respect to subspaces, see Theorem 2.16 of \cite{E}), $C$ is non-countably compact, and thus is non-pseudocompact. (As pointed out, under $\mathbf{MA}+\neg\mathbf{CH}$, any non-compact perfectly normal space is non-countably compact. Also, countable compactness coincides with pseudocompactness in the class of normal spaces, see Theorems 3.10.20 and 3.10.21 of \cite{E}.) Therefore, there is a nonempty zero-set $Z_C$ in $\beta C$ which misses $C$. Observe that $\beta C=\mathrm{cl}_{\beta X}C$, as $C$, being a component of the locally connected space $X$, is both open and closed in $X$. Let $f_C:\mathrm{cl}_{\beta X}C\rightarrow[0,1]$ be a continuous mapping such that $Z(f_C)=Z_C$. Let
\[T=\bigcup_{\mbox{$C$ is a component of $X$}}\mathrm{cl}_{\beta X}C\]
and let
\[f=\sum_{\mbox{$C$ is a component of $X$}}f_C:T\longrightarrow[0,1].\]
Note that distinct components of $X$, being disjoint open and closed subspaces of $X$, have disjoint closures in $\beta X$. Therefore, the mappings $f_C$, where $C$ is a component of $X$, have disjoint domains, and thus $f$ is a well defined mapping. Moreover, $f$ is continuous, as its restriction $f|_{\mathrm{cl}_{\beta X}C}$ ($=f_C$) is continuous and $\mathrm{cl}_{\beta X}C$ is open in its domain for every component $C$ of $X$. Note that $\beta T=\beta X$, as $T$ contains $X$ and is contained in $\beta X$. The mapping $f$ may therefore be extended to a continuous mapping $F:\beta X\rightarrow[0,1]$. Let $Z=Z(F)$. Then, $Z$ is closed in $\beta X$, and is thus compact, and $Z$ is a  $G_\delta$-set in $\beta X$, as $Z$ is a zero-set in $\beta X$. Also, $Z\cap X=\emptyset$, as $Z\cap C=Z(f_C)\cap C=\emptyset$ for every component $C$ of $X$. It is clear that $Z\cap\mathrm{cl}_{\beta X}C=Z(f_C)$ for every component $C$ of $X$, in particular, $Z\cap\mathrm{cl}_{\beta X}C$ is non-empty.

This completes the proof of the claim.
\end{proofoftheclaim}

Assume that $Z$ is a non-empty compact subspace of $\beta X$ which is disjoint from $X$ and intersects the closure in $\beta X$ of every component of $X$ in the cases when $i=3\frac{1}{2},4,5$. In the case when $i=6$ (that we also assume $\mathbf{MA}+\neg\mathbf{CH}$) assume further that $Z$ is a $G_\delta$-set in $\beta X$. Let $Q$ be the quotient space of $\beta X$ obtained by contracting $Z$ to a point $z$, and let $q:\beta X\rightarrow Q$ be the corresponding quotient mapping. Observe that a space obtained from a normal space by contracting a closed subspace to a point is always normal. Thus, in particular, $Q$ is normal. It can be easily checked that $Q$ contains $X$ (in its original topology) as a subspace and $X$ is dense in $Q$. Let $Y=X\cup\{z\}$, considered as a subspace of $Q$. Then, $Y$ is a completely regular space, as is a subspace of $Q$, and $Q$ is completely regular (and complete regularity is hereditary with respect to subspaces, see Theorem 2.1.6 of \cite{E}), and $Y$ contains $X$ as a dense subspace, as $Q$ does. We check that $Y$ is a $T_i$-space and is connected.

\begin{xclaim}
$Y$ is a $T_i$-space.
\end{xclaim}

\begin{proofoftheclaim}
Note that it already follows from the above observation that $Y$ is a completely regular space. In particular, the claim holds for the case when $i=3\frac{1}{2}$.

Let $i=4$. Let $G$ and $H$ be disjoint closed subspaces of $Y$. Suppose that neither $G$ nor $H$ has $z$. Then $G$ and $H$ are disjoint closed subspaces of $X$ and therefore, by normality of $X$, are contained in disjoint open subspaces of $X$ (which are also open subspaces in $Y$). Next, suppose that either $G$ or $H$, say $G$, has $z$. Since $Y$ is regular, there are disjoint open subspaces $U$ and $V$ of $Y$ containing $z$ and $H$, respectively. On the other hand, $G\cap X$ and $H$ are disjoint closed subspaces of $X$ and therefore, by normality of $X$, there are disjoint open subspaces $U'$ and $V'$ of $X$ containing $G\cap X$ and $H$, respectively. Now, $U\cup U'$ and $V\cap V'$ are disjoint open subspaces of $Y$ containing $G$ and $H$, respectively.

Let $i=5$. Note that by the above argument $Y$ is normal. A proof analogous to the one we have given for Theorem \ref{HFD} (case $i=5$) now proves that $Y$ is hereditarily normal.

Let $i=6$ (and assume $\mathbf{MA}+\neg\mathbf{CH}$). Note that by the above argument $Y$ is normal. We check that any closed subspace of $Y$ is a $G_\delta$-set in $Y$. Let $G$ be a closed subspaces of $Y$. Suppose that $G$ does not have $z$. Then $G$ is a closed subspace of $X$, and since $X$ is perfectly normal, is a $G_\delta$-set in $X$ and therefore is a $G_\delta$-set in $Y$. Now, suppose that $G$ has $z$. Recall that by our assumption $Z$ is a $G_\delta$-set in $\beta X$ and so, its image $\{z\}$ under $q$ is a $G_\delta$-set in $Q$ and thus is a $G_\delta$-set in $Y$. Also, $G\cap X$ is a closed subspace of $X$ (and since $X$ is perfectly normal) is a $G_\delta$-set in $X$ and therefore is a $G_\delta$-set in $Y$. Let
\[\{z\}=\bigcap_{n=1}^\infty U_n\quad\mbox{and}\quad G\cap X=\bigcap_{n=1}^\infty V_n, \]
where $U_n$ and $V_n$ are open subspaces of $Y$ such that $U_n\supseteq U_{n+1}$ and $V_n\supseteq V_{n+1}$ for all positive integers $n$. Now
\[G=\bigcap_{n=1}^\infty (U_n\cup V_n),\]
and therefore $G$ is a $G_\delta$-set in $Y$.

This completes the proof of the claim.
\end{proofoftheclaim}

\begin{xclaim}
$Y$ is a connected space.
\end{xclaim}

\begin{proofoftheclaim}
Let $C$ be a component of $X$ and let $Z_C=Z\cap\mathrm{cl}_{\beta X}C$. Then $C\cup Z_C$ is connected, as lies between the connected space $C$ and its closure $\mathrm{cl}_{\beta X}C$ in $\beta X$. In particular, this implies that $q(C\cup Z_C)$, being the continuous image of a connected space, is connected. Observe that $Z_C$ is non-empty by the choice of $Z$ and $q(Z_C)=z$, as $Z_C$ is contained in $Z$ (which is mapped by $q$ to the point $z$). Therefore $q(C\cup Z_C)=C\cup\{z\}$ by the definition of $q$, as $q$ keeps the points of $X$ fixed. But
\[\bigcup_{\mbox{$C$ is a component of $X$}}\big(C\cup\{z\}\big)=X\cup\{z\}=Y,\]
and therefore $Y$, being the union of connected subspaces with the point $z$ in common, is connected.

This completes the proof of the claim.
\end{proofoftheclaim}

This shows that $Y$ is a $T_i$ one-point connectification of $X$, and completes the proof of the theorem.
\end{proof*}

\begin{xremark}
Observe that we have not used the full strength of our local connectedness assumption in Theorem \ref{PHRS}. Indeed, Theorem \ref{PHRS} (and our other theorems) holds true if we replace the local connectedness assumption by the weaker assumption that all components are open, or, equivalently, every point of the space has a connected open neighborhood (a perhaps more natural counterpart of local compactness here in our context).
\end{xremark}

In our next theorem (which is an extension of Theorem \ref{HFD}), for a topological property $\mathscr{P}$, we provide necessary and sufficient conditions for a locally connected $T_i$-space with $\mathscr{P}$ to have a $T_i$ one-point connectification with $\mathscr{P}$, where $i=3\frac{1}{2},4,5$. (There is an analogous theorem in the case when $i=2$ but requires a separate proof; see Theorem \ref{JHJH}.) The conditions on $\mathscr{P}$ are mild and are satisfied by most of the so called ``covering properties'' (i.e., topological properties defined by means of existence of certain kinds of open subcovers or refinements of open covers of certain types). This includes the cases when $\mathscr{P}$ is the Lindel\"{o}f property, paracompactness, metacompactness, subparacompactness, the para-Lindel\"{o}f property, the $\sigma$-para-Lindel\"{o}f property, weak $\theta$-refinability, $\theta$-refinability (or submetacompactness), weak $\delta\theta$-refinability, and $\delta\theta$-refinability (or the submeta-Lindel\"{o}f property). (See \cite{Bu}, \cite{Steph} and \cite{Va} for definitions and proofs.)

For a topological property $\mathscr{P}$, we say that $\mathscr{P}$ is
\begin{enumerate}
  \item \textit{closed hereditary} if every closed subspace of a space with $\mathscr{P}$ also has $\mathscr{P}$.
  \item \textit{preserved under finite closed sums} if any space which is a finite union of closed subspaces each with $\mathscr{P}$ also has $\mathscr{P}$.
  \item \textit{co-local} if a space $X$ has $\mathscr{P}$ provided that it has a point $p$ and an open base $\mathfrak{B}$ at $p$ such that $X\setminus B$ has $\mathscr{P}$ for all $B$ in $\mathfrak{B}$.
\end{enumerate}

Condition (3) in the above definition seems to be first introduced by S. Mr\'{o}wka in \cite{Mr}, where it was referred to as ``condition $(\mathrm{W})$.''

In the proof of the following lemma we will make use of some of the ideas we have already used in the proof of Theorem 2.9 in \cite{Ko}.

\begin{lemma}\label{HHG}
Let $i=3\frac{1}{2},4,5$. Let $\mathscr{P}$ be a closed hereditary co-local topological property preserved under finite closed sums. A locally connected $T_i$-space has a $T_i$ one-point connectification with $\mathscr{P}$ if its components are all non-compact with $\mathscr{P}$.
\end{lemma}

\begin{proof*}
Let $X$ be a locally connected $T_i$-space whose components are all non-compact with $\mathscr{P}$. Let $Z$ be the subspace of $\beta X$ as defined in the first claim in the proof of Theorem \ref{PHRS}. We also let the quotient mapping $q:\beta X\rightarrow Q$ and the $T_i$ one-point connectification $Y=X\cup\{z\}$ of $X$ (which is a subspace of $Q$) be as in the proof of Theorem \ref{PHRS}. To show that $Y$ has $\mathscr{P}$, since $\mathscr{P}$ is a co-local topological property, it suffices to check that $Y\setminus U$ has $\mathscr{P}$ for every open neighborhood $U$ of $z$ in $Y$. So, let $U$ be an open neighborhood of $z$ in $Y$. Let $V$ be an open subspace of $Q$ such that $U=V\cap Y$. Note that $V$ has $z$. Now, $q^{-1}(V)$ is an open subspace of $\beta X$ which contains $q^{-1}(z)$ ($=Z$). Using the definition of $Z$, we have
\[\beta X\setminus q^{-1}(V)\subseteq\beta X\setminus q^{-1}(z)=\beta X\setminus Z\subseteq\bigcup_{\mbox{$C$ is a component of $X$}}\mathrm{cl}_{\beta X}C.\]
Note that $\mathrm{cl}_{\beta X}C$ is open in $\beta X$ for every component $C$ of $X$, as components of $X$ are both open and closed in $X$ (since $X$ is locally connected). On the other hand $\beta X\setminus q^{-1}(V)$ is compact, as is closed in $\beta X$. Therefore
\[\beta X\setminus q^{-1}(V)\subseteq\mathrm{cl}_{\beta X}C_1\cup\cdots\cup\mathrm{cl}_{\beta X}C_n\]
for some components $C_1,\ldots,C_n$ of $X$. We intersect the two sides of the above relation with $X$ to obtain
\[Y\setminus U=X\setminus V\subseteq C_1\cup\cdots\cup C_n.\]
Note that $C_1,\ldots,C_n$ are closed subspaces of $X$ and they all have $\mathscr{P}$ by our assumption. Therefore, their union has $\mathscr{P}$, as $\mathscr{P}$ is preserved under finite closed sums, and so does the closed subspace $Y\setminus U$ of their union, as $\mathscr{P}$ is closed hereditary.
\end{proof*}

Our theorem now follows from Lemma \ref{HHG}.

\begin{theorem}\label{JHG}
Let $i=3\frac{1}{2},4,5$. Let $\mathscr{P}$ be a closed hereditary co-local topological property preserved under finite closed sums. A locally connected $T_i$-space with $\mathscr{P}$ has a $T_i$ one-point connectification with $\mathscr{P}$ if and only if it has no compact component.
\end{theorem}

\begin{proof*}
Sufficiency follows from Lemma \ref{HHG} and the simple observation that components of a space $X$ with $\mathscr{P}$, being closed in $X$, have $\mathscr{P}$. Necessity follows from Theorem \ref{PHRS}.
\end{proof*}

In Theorem \ref{JHG}, the requirements on the topological property $\mathscr{P}$ are rather mild, however, there are still topological properties not satisfying the requirements of Theorem \ref{JHG} for which the conclusion in Theorem \ref{JHG} holds. The most important examples of such topological properties are metrizability, realcompactness and pseudocompactness. Metrizability is closed hereditary and is preserved under finite closed sums but is not co-local, realcompactness is closed hereditary and is co-local but is not preserved under finite closed sums, and pseudocompactness is preserved under finite closed sums and is co-local but is not closed hereditary. (The fact that metrizability is preserved under finite closed sums but is not co-local is shown in Example \ref{KHD}. The fact that realcompactness is closed hereditary is well known; see e.g. Theorem 3.11.4 of \cite{E}. That realcompactness is co-local follows from Corollary 2.6 of \cite{MRW} and Theorems 8.10, 8.11 and 8.16 of \cite{GJ}. Note that in \cite{MRW} the authors call a co-local topological property a topological property which ``satisfy Mr\'{o}wka's condition (W).'' That realcompactness is not preserved under finite closed sums in the realm of completely regular spaces is known. In \cite{Mr1}, with a correction in \cite{Mr2}, S. Mr\'{o}wka describes a non-realcompact completely regular space which is the union of two closed realcompact subspaces; a simpler such a space is latter discovered by A. Mysior in \cite{My}. The facts that pseudocompactness is preserved under finite closed sums and is co-local follow directly from the definitions. That pseudocompactness is not closed hereditary follows from Example 3.10.29 of \cite{E}.)

In the next theorem we show that Theorem \ref{JHG} holds true if $\mathscr{P}$ is metrizability (indeed, a stronger version of Theorem \ref{JHG} holds, as metrizable spaces are perfectly normal), i.e., a locally connected metrizable space has a metrizable one-point connectification if and only if it has no compact component. The fact that Theorem \ref{JHG} (indeed, a stronger version) holds if $\mathscr{P}$ is realcompactness or pseudocompactness will be consecutively shown in Theorems \ref{PIKJ} and \ref{JHOG}.

Observe that pseudocompactness coincides with countable compactness in the class of normal spaces and countable compactness coincides with compactness in the class of metrizable space. (See Theorems 3.10.20, 3.10.21 and 4.1.17 of \cite{E}.) Therefore, for a non-compact metrizable space $X$ there exists a non-empty zero-set of $\beta X$ which is disjoint from $X$.

Recall that by the Nagata--Smirnov metrization theorem, a regular space is metrizable if and only if it has a countably locally finite base (i.e., a base which is a countable union of locally finite collections).

\begin{theorem}\label{PAGA}
A locally connected metrizable space has a metrizable one-point connectification if and only if it has no compact component.
\end{theorem}

\begin{proof*}
As remarked previously in the proof of Theorem \ref{PHRS} a locally connected space which has a Hausdorff one-point connectification has no compact component. We therefore suffice to prove the sufficiency part of the theorem.

Let $X$ be a locally connected metrizable space with no compact component. The proof we have given for Theorem \ref{PHRS} (case $i=6$) now applies to show that there is a non-empty zero-set $Z$ of $\beta X$ which is disjoint from $X$ and intersects the closure in $\beta X$ of every component of $X$. (We need to use the fact that for a non-compact metrizable space $C$ there exists a non-empty zero-set of $\beta C$ which is disjoint from $C$.) Let $F:\beta X\rightarrow[0,1]$ be a continuous mapping such that $Z=Z(F)$. Let the quotient mapping $q:\beta X\rightarrow Q$ and the completely regular one-point connectification $Y=X\cup\{z\}$ of $X$ (which is a subspace of $Q$) be as in the proof of Theorem \ref{PHRS}. To conclude the proof we need to check that $Y$ is metrizable.

By the Nagata--Smirnov metrization theorem, it suffices to check that $Y$ has a countably locally finite base. Note that the space $X$, being metrizable, has a countably locally finite base $\mathfrak{B}$, again by the Nagata--Smirnov metrization theorem. Let
\[\mathfrak{B}=\bigcup_{n=1}^\infty\mathfrak{B}_n,\]
where $\mathfrak{B}_n$ is locally finite for all positive integers $n$. Let
\[\mathfrak{C}=\bigcup_{n=1}^\infty\big\{Y\cap q\big(F^{-1}\big([0,1/n)\big)\big)\big\}\cup\bigcup_{n=1}^\infty\bigcup_{m=1}^\infty\big\{B\setminus q\big(F^{-1}\big([0,1/m]\big)\big):B\in\mathfrak{B}_n\big\}.\]
Then, $\mathfrak{C}$ is a base for $Y$ with the required properties, as one can check.
\end{proof*}

The next example shows that metrizability is not a co-local topological property, though is a closed hereditary  topological property (indeed, is hereditary with respect to arbitrary subspaces) which is preserved under finite closed sums. This in particular shows that Theorem \ref{PAGA} is not deducible from Theorem \ref{JHG}.

\begin{example}\label{KHD}
Let $\mathbb{N}$ be the set of all natural numbers equipped with the discrete topology. Choose some $p$ in $\beta\mathbb{N}\setminus\mathbb{N}$ and let $X=\mathbb{N}\cup\{p\}$, considered as a subspace of $\beta\mathbb{N}$. It is clear that $X\setminus U$, being a subspace of $\mathbb{N}$, is metrizable for every open neighborhood $U$ of $p$ in $X$. We show that $X$ is not metrizable. This will show that metrizability is not a co-local topological property. Suppose to the contrary that $X$ is metrizable. Let $\{U_1,U_2,\ldots\}$ be a countable base for $X$ at $p$. Without any loss of generality we may assume that $U_1\supseteq U_2\supseteq\ldots$ where all inclusions are proper. Define a sequence $k_1,k_2,\ldots$ such that $k_n$ is in $U_n\setminus U_{n+1}$ for all positive integers $n$. Let
\[O=\{k_1,k_3,\ldots\}\quad\mbox{and}\quad E=\{k_2,k_4,\ldots\}.\]
Observe that for all positive integers $n$ the set $U_n$ contains both $O$ and $E$ except for a finite number of elements, and thus, in particular, $\mathrm{cl}_{\beta\mathbb{N}}O\subseteq\mathbb{N}\cup\mathrm{cl}_{\beta\mathbb{N}} U_n$ and $\mathrm{cl}_{\beta\mathbb{N}}E\subseteq\mathbb{N}\cup\mathrm{cl}_{\beta\mathbb{N}} U_n$. Therefore,
\begin{equation}\label{JGH}
\mathrm{cl}_{\beta\mathbb{N}}O\subseteq\mathbb{N}\cup P\quad\mbox{and}\quad \mathrm{cl}_{\beta\mathbb{N}}E\subseteq\mathbb{N}\cup P,
\end{equation}
where $P=\bigcap_{n=1}^\infty\mathrm{cl}_{\beta\mathbb{N}} U_n$. We check that $P=\{p\}$. Clearly $p$ is contained in $P$, as $U_n$'s are all neighborhoods of $p$. Let $q$ be an element of $\beta\mathbb{N}$ distinct from $p$. Let $W$ be an open neighborhood of $p$ in $\beta\mathbb{N}$ whose closure $\mathrm{cl}_{\beta\mathbb{N}}W$ does not contain $q$. Then $W\cap X$ is an open neighborhood of $p$ in $X$, and therefore contains $U_m$ for some positive integer $m$. But then $q$ is not contained in $\mathrm{cl}_{\beta\mathbb{N}}U_m$, as the latter is contained in $\mathrm{cl}_{\beta\mathbb{N}}W$. Thus $q$ is not in $P$. This shows that $P=\{p\}$. Observe that $\mathrm{cl}_{\beta\mathbb{N}}O$ and $\mathrm{cl}_{\beta\mathbb{N}}E$, being compact infinite spaces, are not contained in $\mathbb{N}$. Therefore, by (\ref{JGH}), they both contain $p$. But this is not possible, as $O$ and $E$, being disjoint open and closed subspaces of $\mathbb{N}$, have disjoint closures in $\beta\mathbb{N}$.

Observe that metrizability is, however, hereditary with respect to closed subspaces (indeed, is hereditary with respect to arbitrary subspaces) and is preserved under finite closed sums. (Note that metrizability is an additive topological property, in the sense that, any direct sum of metrizable spaces is metrizable. Now, suppose that a space $X$ is a finite union of closed metrizable subspaces $X_1,\ldots,X_n$. Then, since the mapping
\[\bigcup_{i=1}^n\mathrm{id}_{X_i}:\bigoplus_{i=1}^n X_i\longrightarrow X\]
is perfect continuous and metrizability is preserved under perfect continuous mappings, the space $X$ is metrizable. See Theorems 3.7.22, 4.2.1, 4.4.15 of \cite{E}.)
\end{example}

In the next two theorems we show that Theorem \ref{JHG} (with the inclusion of the additional case $i=6$ under $\mathbf{MA}+\neg\mathbf{CH}$) holds if $\mathscr{P}$ is either realcompactness or pseudocompactness.

\begin{theorem}\label{PIKJ}
Let $i=3\frac{1}{2},4,5$ or $i=6$ with $\mathbf{MA}+\neg\mathbf{CH}$. A locally connected realcompact $T_i$-space has a realcompact $T_i$ one-point connectification if and only if it has no compact component.
\end{theorem}

\begin{proof*}
Let $X$ be a locally connected realcompact $T_i$-space. By Theorem \ref{PHRS} it is clear that if $X$ has a realcompact $T_i$ one-point connectification then it has no compact component. Suppose that $X$ has no compact component. Then, by Theorem \ref{PHRS}, $X$ has a $T_i$ one-point connectification $Y$. Note that by Theorem 8.16 of \cite{GJ} any completely regular space which is the union of a realcompact subspace and a compact subspace is realcompact. Therefore, $Y$ (being the union of $X$ and the singleton $Y\setminus X$) is also realcompact.
\end{proof*}

\begin{theorem}\label{JHOG}
Let $i=3\frac{1}{2},4,5$ or $i=6$ with $\mathbf{MA}+\neg\mathbf{CH}$. A locally connected pseudocompact $T_i$-space has a pseudocompact $T_i$ one-point connectification if and only if it has no compact component.
\end{theorem}

\begin{proof*}
Let $X$ be a locally connected pseudocompact $T_i$-space. By Theorem \ref{PHRS} it is clear that if $X$ has a pseudocompact $T_i$ one-point connectification then it has no compact component. Suppose that $X$ has no compact component. By Theorem \ref{PHRS}, $X$ has a $T_i$ one-point connectification $Y$. Note that $Y$ is also pseudocompact, as it contains the pseudocompact space $X$ as a dense subspace.
\end{proof*}

In the remainder of this section we prove counterparts for Theorems \ref{PHRS} and \ref{JHG} in the case $i=2$.

Recall that a space $X$ is called \textit{almost compact} if for every open cover $\mathscr{U}$ of $X$ there is a finite subcollection $\mathscr{V}$ of $\mathscr{U}$ such that $X=\mathrm{cl}_X\bigcup\mathscr{V}$.

The following theorem is a counterparts for Theorems \ref{PHRS} in the case $i=2$.

\begin{theorem}\label{JLJH}
Let $i=2$. A locally connected $T_i$-space has a $T_i$ one-point connectification if and only if it has no almost compact component.
\end{theorem}

\begin{proof*}
Let $X$ be a locally connected Hausdorff space.

Suppose that $X$ has a Hausdorff one-point connectification $Y$. Suppose to the contrary that $X$ has an almost compact component $C$. Note that $C$ is open in $Y$, as is open in $X$ (since $X$ is locally connected) and $X$ is open in $Y$. We show that $C$ is also closed in $Y$. Let $y$ be in $Y\setminus C$. For every $x$ in $C$ let $U_x$ and $V_x$ be disjoint open neighborhoods of $y$ and $x$ in $Y$, respectively. Let
\[\mathscr{U}=\{C\cap V_x:x\in C\}.\]
Then $\mathscr{U}$ is an open cover of $C$ and therefore the closure in $C$ (and thus the closure in $Y$) of the union of a finite number of elements from $\mathscr{U}$ covers $C$. Let $x_1,\ldots,x_n$ in $C$ be such that $C\subseteq\bigcup_{i=1}^n\mathrm{cl}_Y V_{x_i}$. Let $U=\bigcap_{i=1}^n U_{x_i}$. Then $U$ is an open neighborhood of $y$ in $Y$ which misses $C$. Thus, $y$ is in $Y\setminus\mathrm{cl}_YC$. Therefore $C$ is closed in $Y$. This contradicts connectedness of $Y$.

Now, suppose that $X$ has no almost compact component. Then, for every component $C$ of $X$ there is a cover $\mathscr{U}(C)$ of $C$ consisting of sets which are open in $C$ (and thus open in $X$, as $C$ is open in $X$) such that $C\neq\mathrm{cl}_X\bigcup\mathscr{U}$ for every finite subcollection $\mathscr{U}$ of $\mathscr{U}(C)$.

We proceed with the construction of a Hausdorff one-point connectification for $X$. Take some $p$ outside $X$ and let $Y=X\cup\{p\}$. Let the topology $\mathscr{T}$ of $Y$ consists of (i) all open subspaces of $X$, together with (ii) all subsets $U$ of $Y$ such that
\begin{enumerate}
  \item $p$ is in $U$;
  \item $U\cap X$ is open in $X$;
  \item for every component $C$ of $X$ there is a finite subcollection $\mathscr{U}$ of $\mathscr{U}(C)$ such that \[C\setminus\mathrm{cl}_X\bigcup\mathscr{U}\subseteq U.\]
\end{enumerate}
We verify that this defines a Hausdorff topology on $Y$.

The empty set is of type (i) and the whole set $Y$ is of type (ii). We check that the intersection of any two elements of $\mathscr{T}$ is in $\mathscr{T}$. It is clear that the intersection of any two elements of type (i) is of type (i), and the intersection of an element of type (i) with an element of type (ii) is of type (i). Let $U$ and $V$ be elements of type (ii). It is clear that $p$ is in $U\cap V$ and $(U\cap V)\cap X$ is open in $X$. For a component $C$ of $X$ there are finite subcollections $\mathscr{U}$ and $\mathscr{V}$ of $\mathscr{U}(C)$ such that
\[C\setminus\mathrm{cl}_X\bigcup\mathscr{U}\subseteq U\quad\mbox{and}\quad C\setminus\mathrm{cl}_X\bigcup\mathscr{V}\subseteq V.\]
Let $\mathscr{W}=\mathscr{U}\cup\mathscr{V}$. Then $\mathscr{W}$ is a finite subcollections of $\mathscr{U}(C)$ and
\[C\setminus\mathrm{cl}_X\bigcup\mathscr{W}=\Big(C\setminus\mathrm{cl}_X\bigcup\mathscr{U}\Big)\cap \Big(C\setminus\mathrm{cl}_X\bigcup\mathscr{V}\Big)\subseteq U\cup V.\]
The fact that the union of any number of elements of $\mathscr{T}$ is in $\mathscr{T}$ may be checked analogously. Therefore $\mathscr{T}$ is a topology on $Y$.

Note that $X$ is a subspace of $Y$, as any open set in $X$ is of type (i) and any open set in $Y$ intersected with $X$ is open in $X$. Also, $X$ is dense in $Y$, as any open neighborhood of $p$ in $Y$ contains a set of the form $C\setminus\mathrm{cl}_X\bigcup\mathscr{U}$, where $C$ is a component of $X$ and $\mathscr{U}$ is a finite subcollection of $\mathscr{U}(C)$, and such a set is necessarily non-empty by the definition of $\mathscr{U}(C)$.

We check that $Y$ is Hausdorff. Let $y$ and $z$ be distinct elements in $Y$. First, suppose that $y$ and $z$ are both in $X$. Then $y$ and $z$ are separated by disjoint sets which are open in $X$, and therefore open in $Y$. Next, suppose that either $y$ or $z$, say $y$, is $p$. Then $z$ is in $X$ and therefore is in some component $C$ of $X$. Let $U$ be an element of the cover $\mathscr{U}(C)$ of $C$ which has $z$. Then
\[(C\setminus\mathrm{cl}_X U)\cup (Y\setminus C)\quad\mbox{and}\quad  U\]
are disjoint open neighborhoods of $y$ and $z$ in $Y$, respectively.

To conclude the proof we check that $Y$ is connected. Suppose otherwise that there is a separation $U$ and $V$ for $Y$. Without any loss of generality we may assume that $p$ is in $U$. Then $U$ intersects any component $C$ of $X$, as it contains a non-empty set of the form $C\setminus\mathrm{cl}_X\bigcup\mathscr{U}$ for a finite subcollection $\mathscr{U}$ of $\mathscr{U}(C)$, and thus contains the whole $C$, as $C$ is connected. This implies that $U=Y$, which is a contradiction.
\end{proof*}

The following lemma is a counterpart for Lemma \ref{HHG} in the case $i=2$.

\begin{lemma}\label{JHHG}
Let $i=2$. Let $\mathscr{P}$ be a closed hereditary co-local topological property preserved under finite closed sums. A locally connected $T_i$-space has a $T_i$ one-point connectification with $\mathscr{P}$ if its components are all non-almost compact with $\mathscr{P}$.
\end{lemma}

\begin{proof*}
Let $X$ be a locally connected Hausdorff space. Let components of $X$ be all non-almost compact and have $\mathscr{P}$. For each component $C$ of $X$ let $\mathscr{U}(C)$ be as defined in the proof of Theorem \ref{JLJH}. Let $Y=X\cup\{p\}$, where $p$ is not in $X$. Let the topology of $Y$ consists of all open subspaces of $X$, together with all subsets $U$ of $Y$ satisfying conditions (1)--(3) in the proof of Theorem \ref{JLJH} and
\begin{enumerate}
  \item[(4)] there are components $C_1,\ldots,C_n$ of $X$ such that $Y\setminus(C_1\cup\cdots\cup C_n)\subseteq U$.
\end{enumerate}
This defines a topology on $Y$ and the space $Y$ endowed with this topology is a Hausdorff one-point connectification of $X$. (The proof is essentially the same as the one we have given for Theorem \ref{JLJH}.) We need to check that $Y$ has $\mathscr{P}$. To do this, since $\mathscr{P}$ is co-local, we suffice to check that $Y\setminus U$ has $\mathscr{P}$ for any open neighborhood $U$ of $p$ in $Y$. Let $U$ be an open neighborhood of $p$ in $Y$. By (4), then $Y\setminus U\subseteq C_1\cup\cdots\cup C_n$ for some components $C_1,\ldots,C_n$ of $X$. Note that $C_i$, for each $i=1,\ldots,n$, being a component of $X$, is closed in $X$, and has $\mathscr{P}$ by our assumption. Thus, the finite union $C_1\cup\cdots\cup C_n$ has $\mathscr{P}$, as $\mathscr{P}$ is preserved under finite closed sums. This implies that $Y\setminus U$ has $\mathscr{P}$, as is closed in $C_1\cup\cdots\cup C_n$ and $\mathscr{P}$ is closed hereditary.
\end{proof*}

The following theorem which is a counterpart for Theorem \ref{JHG} in the case $i=2$ now follows almost immediately.

\begin{theorem}\label{JHJH}
Let $i=2$. Let $\mathscr{P}$ be a closed hereditary co-local topological property preserved under finite closed sums. A locally connected $T_i$-space with $\mathscr{P}$ has a $T_i$ one-point connectification with $\mathscr{P}$ if and only if it has no almost compact component.
\end{theorem}

\begin{proof*}
Sufficiency follows from Lemma \ref{JHHG} and the simple observation that components of a space $X$ with $\mathscr{P}$, being closed in $X$, have $\mathscr{P}$. Necessity follows from Theorem \ref{JLJH}.
\end{proof*}

In Theorems \ref{PHRS} and \ref{JLJH} we have excluded the cases when $i=3$. The existing duality between Theorem \ref{HFD} on the one hand and Theorems \ref{PHRS} and \ref{JLJH} on the other hand, however, suggests to conjecture that Theorem \ref{PHRS} or \ref{JLJH} holds in this case as well. We formulate this as a conjecture in the following for possible future reference.

\begin{xconjecture}\label{JHGH}
Let $i=3$. A locally connected $T_i$-space has a $T_i$ one-point connectification if and only if it has no compact (or almost compact) component.
\end{xconjecture}

\begin{xremark}
Observe that in the above conjecture we have not included the cases when $i=0,1$. These cases are of no interest. (This agrees with the analogous case of compactifications which are of particular interest when they are assumed to be Hausdorff.) Indeed, any non-empty $T_0$-space has a $T_0$ one-point connectification (and, of course, any space with a $T_0$ one-point connectification is $T_0$). To see this, let $X$ be a $T_0$-space and let $Y$ be the set obtained by adding one extra point to $X$. Define a topology on $Y$ by adding to the topology of $X$ only one new open subspace, namely, $Y$ itself. The resulting space $Y$ is a $T_0$ one-point connectification of $X$.

Also, a $T_1$-space has a $T_1$ one-point connectification if and only if it is infinite with no isolated point. To see this, let $X$ be an infinite $T_1$-space with no isolated point. Let $Y=X\cup\{p\}$, where $p$ is not in $X$. Let the topology of $Y$ consists of open subspaces of $X$ together with sets of the form $Y\setminus F$, where $F$ is a finite subset of $X$. The fact that this defines a $T_1$ topology on $Y$ is easy to check. It is also easy to check that $Y$ contains $X$ as a subspace, and $X$ is dense in $Y$. (The latter is because we are assuming that $X$ is infinite.) To check that $Y$ is connected, note that the existence of a proper simultaneously open and closed neighborhood of $p$ in $Y$ implies the existence of a non-empty finite open subspace $F$ of $X$. Let $F=\{x_1,\ldots,x_n\}$, where $x_1,\ldots,x_n$ are distinct. Then $\{x_1\}=F\cap(Y\setminus\{x_2,\ldots,x_n\})$, and thus $\{x_1\}$ is open in $X$. This contradicts the fact that $X$ has no isolated point. Therefore $Y$ is connected. For the converse, let $X$ be a $T_1$-space, which is either infinite with an isolated point or finite. Then, in either cases, there is an isolated point $x$ in $X$. Suppose to the contrary that $X$ has a $T_1$ one-point connectification $Y$. Then $\{x\}$ is closed in $Y$, as $Y$ is $T_1$, and $\{x\}$ is open in $Y$, as $\{x\}$ is open in $X$ and $X$ is open in $Y$ (since $Y\setminus X$ is a singleton and is therefore closed in $Y$). This contradicts connectedness of $Y$.
\end{xremark}

\section{The collection of all one-point connectifications of a space as a lattice}

In the previous section we have considered the problem of existence of a one-point connectification by finding a necessary and sufficient condition for a locally connected space to have a one-point connectification. This has been done in the classes of $T_i$-spaces, where $i=3\frac{1}{2},4,5$. (We have extended our theorem to the case when $i=6$ by assuming some set-theoretic assumption, and to the case $i=2$ by slightly modifying its statement. We have also extended our theorem further by proving it in the realm of paracompact, metrizable, realcompact and pseudocompact spaces.) The idea was to propose a duality between the Alexandroff one-point compactification theorem and our one-point connectification theorems. It turns out that despite the existing duality, there are indeed differences. The one-point compactification of a locally compact space is necessarily unique, however, the one-point connectification of a locally connected space may not. The purpose of this section is to study the set of all one-point connectifications of a locally connected locally compact space (naturally partially ordered) in the realm of $T_i$-spaces ($i=3\frac{1}{2},4,5$). We study the properties of this partially ordered set, and for a locally connected locally compact $T_i$-space we prove that this set is a compact conditionally complete lattice whose order-structure determines the topology of all Stone--\v{C}ech remainders of components of the space.

We begin with the following lemma which is motivated by (the proof of) Theorem \ref{PHRS}. The proof is essentially a modification of the proofs we have given in \cite{Ko3} (Theorem 2.7) and in \cite{Ko} (Theorem 2.5). We sketch the proof here for the sake of completeness.

Two one-point connectifications $Y$ and $Z$ of a space $X$ are said to be \textit{equivalent} if there is a homeomorphism between $Y$ and $Z$ which fixes the points of $X$. This indeed defines an equivalence relation on the class of all one-point connectifications of the space $X$. We identify equivalence classes with individuals when no confusion arises.

A mapping $f:P\rightarrow Q$ between two partially ordered sets $P$ and $Q$ is called an  \textit{anti-order-homomorphism} if $f(b)\leq f(a)$ whenever $a\leq b$ and $a$ and $b$ are in $P$. The mapping $f$ is an \textit{anti-order-isomorphism} if it is bijective, and $f$ and $f^{-1}$ are both anti-order-homomorphisms.

\begin{notationandlemma}\label{YGFD}
Let $i=3\frac{1}{2},4,5$. Let $X$ be a locally connected $T_i$-space. Let
\[\mathscr{Y}_i(X)=\{Y:\mbox{$Y$ is a $T_i$ one-point connectification of $X$}\},\]
where for two elements $Y_1$ and $Y_2$ of $\mathscr{Y}_i(X)$ we let $Y_1\leq Y_2$ if there is a continuous mapping $f:Y_2\rightarrow Y_1$ which fixes the points of $X$. Let
\[\mathscr{Z}(X)=\{Z\subseteq\beta X\setminus X:\mbox{$Z$ is compact, $Z\cap\mathrm{cl}_{\beta X}C\neq\emptyset$ for all components $C$ of $X$}\},\]
partially ordered with inclusion. For a $Z$ in $\mathscr{Z}(X)$ let $Q$ be the quotient space of $\beta X$ obtained by contracting $Z$ to a point $z$ and let
\[\zeta_X(Z)=X\cup\{z\},\]
considered as a subspace of $Q$. Then $\zeta_X(Z)$ is in $\mathscr{Y}_i(X)$, and the mapping
\[\zeta_X:\big(\mathscr{Z}(X),\subseteq\big)\longrightarrow\big(\mathscr{Y}_i(X),\leq\big)\]
is an anti-order-isomorphism.
\end{notationandlemma}

\begin{proof*}
It is clear that $\leq$ defines a partial order on the set $\mathscr{Y}_i(X)$. It is also clear from the construction of the one-point connectification given in the proof of Theorem \ref{PHRS} that $\zeta_X$ indeed defines a mapping. Note that for a $Y=X\cup\{z\}$ in $\mathscr{Y}_i(X)$, the space $\beta Y$ coincides with the quotient space of $\beta X$ obtained by contracting $\phi^{-1}(z)$ to $z$, where $\phi:\beta X\rightarrow\beta Y$ is the continuous extension of the identity mapping on $X$. We check that $\phi^{-1}(z)$ is in $\mathscr{Z}(X)$. Note that if $\phi^{-1}(z)$ misses $\mathrm{cl}_{\beta X}C$ for some component $C$ of $X$, then $z$ is not in $\phi(\mathrm{cl}_{\beta X}C)$, and therefore $z$ is not in $\mathrm{cl}_YC$, as $\mathrm{cl}_{\beta Y}C\subseteq\phi(\mathrm{cl}_{\beta X}C)$ (since the latter is compact and contains $C$). Thus $C$ is closed in $Y$, as $C$ is closed in $X$. But $C$ is also open in $Y$, as $C$ is open in $X$ (since $X$ is locally connected) and $X$ is open $Y$. This contradicts connectedness of $Y$. Therefore $\zeta_X$ is surjective. The fact that $\zeta_X$ is injective, and $\zeta_X$ and $\zeta_X^{-1}$ are anti-order-homomorphisms follow from standard arguments.
\end{proof*}

A \textit{conditionally complete lattice} is a partially ordered set $\mathscr{L}$ in which every non-empty bounded subset $\mathscr{H}$ has the least upper bound $\bigvee\mathscr{H}$ and the greatest lower bound $\bigwedge\mathscr{H}$. A conditionally complete lattice $\mathscr{L}$ is \textit{compact} if for every subset $\mathscr{H}$ of $\mathscr{L}$ for which $\bigvee\mathscr{H}$ does not exist there is a finite subset $\mathscr{H'}$ of $\mathscr{H}$ for which $\bigvee\mathscr{H'}$ also does not exist.

\begin{theorem}\label{KDSG}
Let $i=3\frac{1}{2},4,5$. Let $X$ be a locally connected locally compact $T_i$-space. The partially ordered set $\mathscr{Y}_i(X)$ of all $T_i$ one-point connectifications of $X$ is a compact conditionally complete lattice.
\end{theorem}

\begin{proof*}
To show that $\mathscr{Y}_i(X)$ is a conditionally complete lattice, by Lemma \ref{YGFD}, it suffices to show that $\mathscr{Z}(X)$ (which is partially ordered with inclusion) is a conditionally complete lattice. Let $\mathscr{H}$ be a bounded subset of $\mathscr{Z}(X)$. Let $Z$ be a lower bound for $\mathscr{H}$ in $\mathscr{Z}(X)$. Then $\bigcap\mathscr{H}$ is a compact subspace of $\beta X\setminus X$ (as each element of $\mathscr{H}$ is) which intersects the closure in $\beta X$ of every component of $X$, as $\bigcap\mathscr{H}$ contains $Z$ and $Z$ does. Therefore $\bigcap\mathscr{H}$ is in $\mathscr{Z}(X)$ which is the greatest lower bound for $\mathscr{H}$ in $\mathscr{Z}(X)$. Note that $\beta X\setminus X$ is closed in $\beta X$, as $X$ is locally compact. In particular, $\beta X\setminus X$, containing $\bigcup\mathscr{H}$, contains its closure $\mathrm{cl}_{\beta X}\bigcup\mathscr{H}$. It is clear that $\mathrm{cl}_{\beta X}\bigcup\mathscr{H}$ intersects the closure in $\beta X$ of every component of $X$ (as each element of $\mathscr{H}$ does) and is therefore in $\mathscr{Z}(X)$. Also, any upper bound of $\mathscr{H}$ in $\mathscr{Z}(X)$ contains the union $\bigcup\mathscr{H}$, and therefore contains the closure $\mathrm{cl}_{\beta X}\bigcup\mathscr{H}$. That is, $\mathrm{cl}_{\beta X}\bigcup\mathscr{H}$ is the least upper bound for $\mathscr{H}$ in $\mathscr{Z}(X)$.

Next, we show that $\mathscr{Y}_i(X)$ is compact. By Lemma \ref{YGFD}, it suffices to check that for every subset $\mathscr{H}$ of $\mathscr{Z}(X)$ for which $\bigcap\mathscr{H}$ is not in $\mathscr{Z}(X)$ there is a finite subset $\mathscr{H'}$ of $\mathscr{H}$ for which $\bigcap\mathscr{H}'$ is not in $\mathscr{Z}(X)$.

So, let $\mathscr{H}$ be a subset of $\mathscr{Z}(X)$ such that $\bigcap\mathscr{H}$ is not in $\mathscr{Z}(X)$. Note that $\bigcap\mathscr{H}$ is a compact subspace of $\beta X\setminus X$. Therefore
\[\mathrm{cl}_{\beta X}C\cap\bigcap\mathscr{H}=\emptyset\]
for some component $C$ of $X$. Compactness now implies that
\[\mathrm{cl}_{\beta X}C\cap\bigcap\mathscr{H}'=\emptyset\]
for some finite subset $\mathscr{H}'$ of $\mathscr{H}$. In particular, then $\bigcap\mathscr{H}'$ is not in $\mathscr{Z}(X)$.
\end{proof*}

The next theorem relates the order-structure of the set of all $T_i$ ($i=3\frac{1}{2},4,5$) one-point connectifications of a locally connected locally compact $T_i$-space to the topology of Stone--\v{C}ech remainders of its components. The first such a result has been obtained by K.D. Magill, Jr. in \cite{Ma} in which the author proved that the order-structure of the set of all compactifications of a locally compact Hausdorff space $X$ determines (and is determined by) the topology of the Stone--\v{C}ech remainder $\beta X\setminus X$. Similar results may be found in \cite{Ko1} (also in \cite{Ko2}) where we have proved that the order-structure of the set of all pseudocompactifications (realcompactifications, respectively) with compact remainder of a locally pseudocompact (locally realcompact, respectively) completely regular space $X$ determines the topology of the subspace $\mathrm{cl}_{\beta X}(\beta X\setminus\upsilon X)$ ($\mathrm{cl}_{\beta X}(\upsilon X\setminus X)$, respectively) of the Stone--\v{C}ech compactification of $X$. We also mention the result obtained by M. Henriksen, L. Janos and R.G. Woods in \cite{HJW} in which the authors proved that the order-structure of the set of all one-point metrizable extensions of a locally compact separable metrizable space $X$ determines (and is determined by) the topology of the Stone--\v{C}ech remainder $\beta X\setminus X$. Other similar results relating the order-structure of certain collections of extensions of a completely regular space $X$ to the topology of certain subspaces of its Stone--\v{C}ech remainder $\beta X\setminus X$ (or remainders in other compactifications) may be found in \cite{D}, \cite{Ko8}, \cite{Ko7}, \cite{MRW}, \cite{Me}, \cite{PW3}, \cite{PW2}, \cite{R}, and \cite{Wo}.

\begin{theorem}\label{QDFK}
Let $i=3\frac{1}{2},4,5$. Let $X$ be a locally connected locally compact $T_i$-space. The order-structure of the set $\mathscr{Y}_i(X)$ of all $T_i$ one-point connectifications of $X$ determines the (faithfully indexed) collection
\[\{\beta C\setminus C:\mbox{$C$ is a component of $X$}\},\]
of all Stone--\v{C}ech remainders of components of $X$, if neither is of cardinality $\leq 2$.
\end{theorem}

\begin{proof*}
Let $X$ and $Y$ be locally connected locally compact $T_i$-spaces such that the partially ordered sets $\mathscr{Y}_i(X)$ and $\mathscr{Y}_i(Y)$ (of all $T_i$ one-point connectifications of $X$ and $Y$, respectively) are order-isomorphic. Also, let the Stone--\v{C}ech remainders of all components of $X$ and $Y$ be of cardinality at least 3. We prove that there is a one-to-one correspondence $\psi$ between the set of all components of $X$ and the set of all components of $Y$, carrying every component $C$ of $X$ to the component $C_\psi$ of $Y$ such that the spaces $\beta C\setminus C$ and $\beta C_\psi\setminus C_\psi$ are homeomorphic. Note that by Lemma \ref{YGFD} the partially ordered sets $\mathscr{Z}(X)$ and $\mathscr{Z}(Y)$ (consisting of all non-empty compact subspaces of $\beta X\setminus X$ and $\beta Y\setminus Y$ which intersect the closure in $\beta X$ and $\beta Y$ of every component of $X$ and $Y$, respectively) are anti-order-isomorphic to $\mathscr{Y}_i(X)$ and $\mathscr{Y}_i(Y)$, respectively, and are therefore themselves order-isomorphic. We let $\phi:\mathscr{Z}(X)\rightarrow\mathscr{Z}(Y)$ denote an order-isomorphism. For simplicity of the notation we let $\mathscr{C}(X)$ and $\mathscr{C}(Y)$ denote the set of all components of $X$ and $Y$, respectively.

We divide the proof into verification of several claims.

Our first step is to identify the minimal elements of $\mathscr{Z}(X)$. (Obviously, here and elsewhere in the proof, because of symmetry, analogous conclusions hold if we consider $X$ instead of $Y$.)

\begin{xclaim}\label{BGFDB}
For an $M$ in $\mathscr{Z}(X)$ the following are equivalent:
\begin{itemize}
\item[(1)] $M$ is minimal in $\mathscr{Z}(X)$.
\item[(2)] We have
\[M=\mathrm{cl}_{\beta X}\bigg(\bigcup_{C\in\mathscr{C}(X)}\{m_C\}\bigg),\]
where $m_C$ is in $\mathrm{cl}_{\beta X}C$ for every component $C$ of $X$.
\end{itemize}
\end{xclaim}

\begin{proofoftheclaim}
(1) \textit{implies} (2). The set $M$, being in $\mathscr{Z}(X)$, intersects the closure in $\beta X$ of every component of $X$. Let $m_C$ be in $M\cap\mathrm{cl}_{\beta X}C$ for every component $C$ of $X$. Let
\[A=\bigcup_{C\in\mathscr{C}(X)}\{m_C\}.\]
Note that $\beta X\setminus X$ contains $A$ and, being closed in $\beta X$ (since $X$ is locally compact), contains its closure $\mathrm{cl}_{\beta X}A$ (denoted by $Z$). Note that $Z$ intersects $\mathrm{cl}_{\beta X}C$ for every component $C$ of $X$, as $m_C$ is in $Z$. Thus $Z$ is in $\mathscr{Z}(X)$. Note that $M$ contains $A$, and therefore contains its closure $Z$ in $\beta X$. Now, minimality of $M$ implies that $Z=M$.

\bigskip

(2) \textit{implies} (1). Let $Z$ be in $\mathscr{Z}(X)$ such that $Z\subseteq M$. Let $D$ be a component of $X$. Note that components of $X$ are open in $X$, as $X$ is locally connected, and thus are open and closed in $X$. In particular,
\[D\quad\text{and}\quad\bigcup_{D\neq C\in\mathscr{C}(X)}C,\]
are disjoint open and closed subspaces of $X$, and therefore, have disjoint closures in $\beta X$. But
\[\mathrm{cl}_{\beta X}\bigg(\bigcup_{D\neq C\in\mathscr{C}(X)}\{m_C\}\bigg)\subseteq\mathrm{cl}_{\beta X}\bigg(\bigcup_{D\neq C\in\mathscr{C}(X)}C\bigg)\]
and therefore
\[\mathrm{cl}_{\beta X}D\cap\mathrm{cl}_{\beta X}\bigg(\bigcup_{D\neq C\in\mathscr{C}(X)}\{m_C\}\bigg)=\emptyset.\]
Observe that
\[M=\{m_D\}\cup\mathrm{cl}_{\beta X}\bigg(\bigcup_{D\neq C\in\mathscr{C}(X)}\{m_C\}\bigg).\]
Thus $M\cap\mathrm{cl}_{\beta X}D=\{m_D\}$. But this implies that $Z\cap\mathrm{cl}_{\beta X}D=\{m_D\}$, as $Z\cap\mathrm{cl}_{\beta X}D$ is non-empty and is contained in $M\cap\mathrm{cl}_{\beta X}D$. Therefore, $Z$ contains the set of all $m_D$ for components $D$ of $X$ and so contains the closure of this set in $\beta X$, namely, $M$. Thus $Z=M$ which proves the minimality of $M$.

This completes the proof of the claim.
\end{proofoftheclaim}

For our purpose, let us call an element $D$ of $\mathscr{Z}(X)$ a \textit{dual atom} if
\[\big|\big\{Z\in\mathscr{Z}(X):Z\subsetneq D\big\}\big|=2.\]

Next, we identify the dual atoms of $\mathscr{Z}(X)$.

\begin{xclaim}\label{GFP}
For a $D$ in $\mathscr{Z}(X)$ the following are equivalent:
\begin{itemize}
\item[(1)] $D$ is a dual atom in $\mathscr{Z}(X)$.
\item[(2)] There is a component $E$ of $X$ and distinct elements $a$ and $b$ in $\mathrm{cl}_{\beta X}E$ such that
\[D=\{a,b\}\cup\mathrm{cl}_{\beta X}\bigg(\bigcup_{E\neq C\in\mathscr{C}(X)}\{d_C\}\bigg),\]
where $d_C$ is in $\mathrm{cl}_{\beta X}C$ for any component $C$ of $X$ distinct from $E$.
\end{itemize}
\end{xclaim}

\begin{proofoftheclaim}
(1) \textit{implies} (2). First, we show that there is a component $E$ of $X$ such that
\begin{equation}\label{JHDS}
|D\cap\mathrm{cl}_{\beta X}E|>1.
\end{equation}
Suppose otherwise. Note that $D\cap\mathrm{cl}_{\beta X}C$ is non-empty, and thus $|D\cap\mathrm{cl}_{\beta X}C|=1$ for every components $C$ of $X$. Denote
\begin{equation}\label{PJGH}
D\cap\mathrm{cl}_{\beta X}C=\{d_C\}
\end{equation}
for every components $C$ of $X$. Let
\[M=\mathrm{cl}_{\beta X}\bigg(\bigcup_{C\in\mathscr{C}(X)}\{d_C\}\bigg).\]
Then $M\subseteq D$, and $M$ is minimal in $\mathscr{Z}(X)$ by the first claim. Therefore $M\neq D$. We consider the following two cases.

\subsubsection*{Case 1.} Suppose that $|D\setminus M|=1$. Then $D=M\cup\{d\}$. Note that $d$ is not in $\mathrm{cl}_{\beta X}C$ for any component $C$ of $X$ by (\ref{PJGH}), as $d$ is not in $M$. Now, let $Z$ be an element of $\mathscr{Z}(X)$ such that $Z\subsetneq D$. Note that $Z\cap\mathrm{cl}_{\beta X}C$ is non-empty, and is contained in $D\cap\mathrm{cl}_{\beta X}C$ for every components $C$ of $X$. Using (\ref{PJGH}), then $d_C$ is in $Z$ for every components $C$ of $X$. Therefore $M\subseteq Z$ which implies that $M=Z$. That is, there is only one element $Z$ in $\mathscr{Z}(X)$ such that $Z\subsetneq D$, namely, $M$.

\subsubsection*{Case 2.} Suppose that $|D\setminus M|>1$. Let $a$ and $b$ be distinct elements in $D\setminus M$. Let $Z$ be either
\[M,\quad M\cup\{a\}\quad\mbox{or}\quad M\cup\{b\}.\]
Then $Z$ is in $\mathscr{Z}(X)$ and $Z\subsetneq D$. That is, there are at least 3 elements $Z$ in $\mathscr{Z}(X)$ such that $Z\subsetneq D$.

\bigskip

\noindent Thus in either cases the number of elements $Z$ of $\mathscr{Z}(X)$ with $Z\subsetneq D$ is not 2. This contradiction shows (\ref{JHDS}) for a component $E$ of $X$.

We show that $E$ is unique with this property. Suppose otherwise. Let $F$ be a component of $X$ distinct from $E$ and such that
\[|D\cap\mathrm{cl}_{\beta X}F|>1.\]
Let
\[\{a,e\}\subseteq D\cap\mathrm{cl}_{\beta X}E\quad\mbox{and}\quad\{b,f\}\subseteq D\cap\mathrm{cl}_{\beta X}F\]
where $a$ and $e$, and $b$ and $f$, are pairs of distinct elements. Let
\[Z=\{u,v\}\cup\mathrm{cl}_{\beta X}\bigg(\bigcup_{E,F\neq C\in\mathscr{C}(X)}\{d_C\}\bigg),\]
where $\{u,v\}$ is either
\[\{a,b\},\quad\{a,f\},\quad\{e,b\},\quad\mbox{or}\quad\{e,f\}.\]
Then the elements $Z$ are in $\mathscr{Z}(X)$, are distinct, and $Z\subsetneq D$. (That the elements $Z$ are distinct follows from comparing $Z\cap\mathrm{cl}_{\beta X}E$ and $Z\cap\mathrm{cl}_{\beta X}F$, using an argument similar to the one given in the proof of the first claim, part (2)$\Rightarrow$(1).) This contradiction shows the uniqueness of $E$.

Next, we show that
\[|D\cap\mathrm{cl}_{\beta X}E|=2.\]
Suppose otherwise. Let $a$, $b$ and $c$ be distinct elements of $D\cap\mathrm{cl}_{\beta X}E$. Let
\[Z=\{z\}\cup\mathrm{cl}_{\beta X}\bigg(\bigcup_{E\neq C\in\mathscr{C}(X)}\{d_C\}\bigg),\]
where $z$ is either $a$, $b$ or $c$. Then such $Z$ are elements of $\mathscr{Z}(X)$, are distinct, and $Z\subsetneq D$. This is a contradiction.

Let
\[D\cap\mathrm{cl}_{\beta X}E=\{a,b\}.\]
Let
\[Z=\{a,b\}\cup\mathrm{cl}_{\beta X}\bigg(\bigcup_{E\neq C\in\mathscr{C}(X)}\{d_C\}\bigg).\]
Then $Z$ is an element of $\mathscr{Z}(X)$ and $Z\subseteq D$. Observe that
\[T=\{t\}\cup\mathrm{cl}_{\beta X}\bigg(\bigcup_{E\neq C\in\mathscr{C}(X)}\{d_C\}\bigg),\]
where $t$ is either $a$ or $b$, are distinct elements of $\mathscr{Z}(X)$ with $T\subsetneq Z$. From this it follows that $D=Z$.

\bigskip

(2) \textit{implies} (1). Let $Z$ be an element of $\mathscr{Z}(X)$ such that $Z\subsetneq D$. Note that $d_C$ is in $Z$ for every component $C$ of $X$ distinct from $E$, as $Z\cap\mathrm{cl}_{\beta X}C$ is non-empty and
\[Z\cap\mathrm{cl}_{\beta X}C\subseteq D\cap\mathrm{cl}_{\beta X}C=\{d_C\}.\]
In particular, then
\[\mathrm{cl}_{\beta X}\bigg(\bigcup_{E\neq C\in\mathscr{C}(X)}\{d_C\}\bigg)\subseteq Z.\]
Also, $Z\cap\mathrm{cl}_{\beta X}E$ is either $\{a\}$, $\{b\}$, or $\{a,b\}$ (where the latter case does not occur because $Z\neq D$), as is non-empty and
\[Z\cap\mathrm{cl}_{\beta X}E\subseteq D\cap\mathrm{cl}_{\beta X}E=\{a,b\}.\]
Therefore
\[Z=\{z\}\cup\mathrm{cl}_{\beta X}\bigg(\bigcup_{E\neq C\in\mathscr{C}(X)}\{d_C\}\bigg),\]
where $z$ is either $a$ or $b$. Obviously, such $Z$ are in $\mathscr{Z}(X)$, are distinct, and $Z\subsetneq D$.

This completes the proof of the claim.
\end{proofoftheclaim}

Next, we define a mapping
\[\psi:\bigcup_{C\in\mathscr{C}(X)}(\mathrm{cl}_{\beta X}C\setminus X)\longrightarrow\bigcup_{D\in\mathscr{C}(Y)}(\mathrm{cl}_{\beta Y}D\setminus Y)\]
and verify that it is a homeomorphism.

Note that components of the space $X$ are pairwise disjoint open and closed subspaces of $X$. Thus, their closures in $\beta X$ are open and closed in $\beta X$ and are pairwise disjoint. Also, $\mathrm{cl}_{\beta X}C=\beta C$, and, in particular
\[\mathrm{cl}_{\beta X}C\setminus X=\beta C\setminus C\]
for every component $C$ of $X$. So the domain of $\psi$ would be the disjoint union (or sum) of $\beta C\setminus C$ for all components $C$ of $X$. The same thing holds for the codomain of $\psi$.

We now proceed with the definition of $\psi$. For every component $C$ of $X$ fix some $m_C$ in $\mathrm{cl}_{\beta X}C\setminus X$. Let
\[M=\mathrm{cl}_{\beta X}\bigg(\bigcup_{C\in\mathscr{C}(X)}\{m_C\}\bigg).\]
Note that $M$ is a minimal element of $\mathscr{Z}(X)$ by our first claim. Let $N=\phi(M)$ (where $\phi:\mathscr{Z}(X)\rightarrow\mathscr{Z}(Y)$ is the order-isomorphism introduced at beginning of the proof).
Then $N$ is a minimal element of $\mathscr{Z}(Y)$ and is therefore of the form
\[N=\mathrm{cl}_{\beta Y}\bigg(\bigcup_{D\in\mathscr{C}(Y)}\{n_D\}\bigg),\]
where $n_D$ is in $\mathrm{cl}_{\beta Y}D\setminus Y$ for every component $D$ of $Y$ by our first claim.

Let $c$ be in $\mathrm{cl}_{\beta X}C\setminus X$ for a component $C$ of $X$. We will define $\psi(c)$. We choose distinct elements $a$ and $b$ in $\mathrm{cl}_{\beta X}C\setminus X$ both distinct from $c$. Note that this is possible by our assumption on cardinalities of the Stone--\v{C}ech remainders of components of $X$. Let
\[A=\{a,c\}\cup\mathrm{cl}_{\beta X}\bigg(\bigcup_{C\neq E\in\mathscr{C}(X)}\{m_E\}\bigg)\quad\mbox{and}\quad B=\{b,c\}\cup\mathrm{cl}_{\beta X}\bigg(\bigcup_{C\neq E\in\mathscr{C}(X)}\{m_E\}\bigg).\]
Then $A$ and $B$ are dual atoms in $\mathscr{Z}(X)$, by our second claim, and therefore, so are their images $\phi(A)$ and $\phi(B)$ under the order-isomorphism $\phi$. By our second claim we have
\[\phi(A)=\{s,t\}\cup\mathrm{cl}_{\beta Y}\bigg(\bigcup_{G\neq F\in\mathscr{C}(Y)}\{k_F\}\bigg)\!\!\!\!\quad\mbox{and}\quad\!\!\!\!\phi(B)=\{u,v\}\cup\mathrm{cl}_{\beta Y}\bigg(\bigcup_{H\neq F\in\mathscr{C}(Y)}\{l_F\}\bigg)\]
where $G$ and $H$ are components of $Y$, $s$ and $t$ are distinct elements in $\mathrm{cl}_{\beta Y}G\setminus Y$, $u$ and $v$ are distinct elements in $\mathrm{cl}_{\beta Y}H\setminus Y$, $k_F$ is in $\mathrm{cl}_{\beta Y}F\setminus Y$ for every component $F$ of $Y$ distinct from $G$, and $l_F$ is in $\mathrm{cl}_{\beta Y}F\setminus Y$ for every component $F$ of $Y$ distinct from $H$.

First, we check that $k_F=n_F$ for every component $F$ of $Y$ distinct from $G$. The same proof applies to show that $l_F=n_F$ for every component $F$ of $Y$ distinct from $H$. Suppose otherwise. Let $I$ be a component of $Y$ distinct from $G$ such that $k_I=n_I$. Observe that
\[A\cup M=\{a,c,m_C\}\cup\mathrm{cl}_{\beta X}\bigg(\bigcup_{C\neq E\in\mathscr{C}(X)}\{m_E\}\bigg).\]
There are at most 3 minimal elements $Z$ in $\mathscr{Z}(X)$ such that $Z\subseteq A\cup M$, namely,
\[Z=\{z\}\cup\mathrm{cl}_{\beta X}\bigg(\bigcup_{C\neq E\in\mathscr{C}(X)}\{m_E\}\bigg),\]
where $z$ is either $a$, $c$, or $m_C$ (which is minimal in $\mathscr{Z}(X)$ by the representation given in our first claim). On the other hand,
\begin{align*}
\phi(A\cup M)&=\phi(A)\cup\phi(M)\\
 &=\{s,t\}\cup\{k_I\}\cup\mathrm{cl}_{\beta Y}\bigg(\bigcup_{G,I\neq F\in\mathscr{C}(Y)}\{k_F\}\bigg)\\
 &\cup\{n_G\}\cup\{n_I\}\cup\mathrm{cl}_{\beta Y}\bigg(\bigcup_{G,I\neq F\in\mathscr{C}(Y)}\{n_F\}\bigg)\\
 &=\{s,t,n_G\}\cup\{k_I,n_I\}\cup\mathrm{cl}_{\beta Y}\bigg(\bigcup_{G,I\neq F\in\mathscr{C}(Y)}\{k_F,n_F\}\bigg).
\end{align*}
Thus, there are at least 4 minimal elements $Z$ in $\mathscr{Z}(Y)$ such that $Z\subseteq\phi(A\cup M)$, namely,
\[Z=\{y\}\cup\{z\}\cup\mathrm{cl}_{\beta Y}\bigg(\bigcup_{G,I\neq F\in\mathscr{C}(Y)}\{z_F\}\bigg),\]
where $y$ is either $s$, $t$, or $n_G$ (from which at least two are distinct), $z$ is either $k_I$ or $n_I$ (which are distinct) and $z_F$ is either $k_F$ or $n_F$ for every component $F$ of $Y$ distinct from $G$ and $I$. This contradiction proves our assertion.

Next, we check that $G=H$. Suppose otherwise. Note that
\[A\cup B=\{a,b,c\}\cup\mathrm{cl}_{\beta X}\bigg(\bigcup_{C\neq E\in\mathscr{C}(X)}\{m_E\}\bigg).\]
In particular, there are 3 minimal elements $Z$ in $\mathscr{Z}(X)$ such that $Z\subseteq A\cup B$, namely,
\[Z=\{z\}\cup\mathrm{cl}_{\beta X}\bigg(\bigcup_{C\neq E\in\mathscr{C}(X)}\{m_E\}\bigg),\]
where $z$ is either $a$, $c$, or $c$. On the other hand,
\begin{align}
\phi(A\cup B)
 &=\phi(A)\cup\phi(B)\nonumber\\
 &=\{s,t\}\cup\{n_H\}\cup\mathrm{cl}_{\beta Y}\bigg(\bigcup_{G,H\neq F\in\mathscr{C}(Y)}\{n_F\}\bigg)\nonumber\\
 &\cup\{n_G\}\cup\{u,v\}\cup\mathrm{cl}_{\beta Y}\bigg(\bigcup_{G,H\neq F\in\mathscr{C}(Y)}\{n_F\}\bigg)\nonumber\\
 &=\{s,t,n_G\}\cup\{u,v,n_H\}\cup\mathrm{cl}_{\beta Y}\bigg(\bigcup_{G,H\neq F\in\mathscr{C}(Y)}\{n_F\}\bigg).
\label{IJGHF}
\end{align}
Thus, there are at least 4 minimal elements $Z$ in $\mathscr{Z}(Y)$ such that $Z\subseteq\phi(A\cup B)$, namely,
\[Z=\{y\}\cup\{z\}\cup\mathrm{cl}_{\beta Y}\bigg(\bigcup_{G,H\neq F\in\mathscr{C}(Y)}\{n_F\}\bigg),\]
where $y$ is either $s$, $t$, or $n_G$ (from which at least two are distinct) and $z$ is either $u$, $v$ or $n_H$ (from which at least two are distinct). This is a contradiction. Therefore $G$ and $H$ are identical. Denote the common value by $D$.

Form (\ref{IJGHF}) it follows that
\[\phi(A\cup B)=\{s,t,u,v\}\cup\mathrm{cl}_{\beta Y}\bigg(\bigcup_{D\neq F\in\mathscr{C}(Y)}\{n_F\}\bigg).\]
As pointed out in above, there are 3 minimal elements $Z$ in $\mathscr{Z}(X)$ with $Z\subseteq A\cup B$. Thus, there will be also 3 minimal elements $Z$ in $\mathscr{Z}(Y)$ such that $Z\subseteq\phi(A\cup B)$. Therefore, from the elements $s$, $t$, $u$, and $v$, only 3 are distinct. But, the pair $s$ and $t$, and the pair $u$ and $v$, both consist of distinct elements. Therefore, the two sets $\{s,t\}$ and $\{u,v\}$ have a single element in common. Define $\psi(c)$ to be the single element in the intersection $\{s,t\}\cap\{u,v\}$. That is,
\[\big\{\psi(c)\big\}=\{s,t\}\cap\{u,v\}.\]

\begin{xclaim}\label{HG}
$\psi$ is a well defined mapping.
\end{xclaim}

\begin{proofoftheclaim}
We need to check that for the element $c$ in $\mathrm{cl}_{\beta X}C\setminus X$, the value $\psi(c)$ is independent of the choices of the elements $a$ and $b$ in $\mathrm{cl}_{\beta X}C\setminus X$. Suppose that $a^*$ and $b^*$ are distinct elements in $\mathrm{cl}_{\beta X}C\setminus X$ both distinct from $c$. Let
\[A^*=\{a^*,c\}\cup\mathrm{cl}_{\beta X}\bigg(\bigcup_{C\neq E\in\mathscr{C}(X)}\{m_E\}\bigg)\]
and
\[B^*=\{b^*,c\}\cup\mathrm{cl}_{\beta X}\bigg(\bigcup_{C\neq E\in\mathscr{C}(X)}\{m_E\}\bigg).\]
Let
\[\phi(A^*)=\{s^*,t^*\}\cup\mathrm{cl}_{\beta Y}\bigg(\bigcup_{D^*\neq F\in\mathscr{C}(Y)}\{n_F\}\bigg)\]
and
\[\phi(B^*)=\{u^*,v^*\}\cup\mathrm{cl}_{\beta Y}\bigg(\bigcup_{D^*\neq F\in\mathscr{C}(Y)}\{n_F\}\bigg),\]
where $D^*$ is a component of $Y$ and $s^*$, $t^*$, $u^*$ and $v^*$ are elements in $\mathrm{cl}_{\beta Y}D^*\setminus Y$.

First, we check that $D=D^*$. Suppose otherwise. Note that
\[A\cup B\cup A^*\cup B^*=\{a,b,a^*,b^*,c\}\cup\mathrm{cl}_{\beta X}\bigg(\bigcup_{C\neq E\in\mathscr{C}(X)}\{m_E\}\bigg).\]
Therefore, there are at most 5 minimal elements $Z$ in $\mathscr{Z}(X)$ such that
\[Z\subseteq A\cup B\cup A^*\cup B^*,\]
namely,
\[Z=\{z\}\cup\mathrm{cl}_{\beta X}\bigg(\bigcup_{C\neq E\in\mathscr{C}(X)}\{m_E\}\bigg)\]
where $z$ is either $a$, $b$, $a^*$, $b^*$ or $c$. On the other hand,
\begin{align*}
&\phi(A)\cup\phi(B)\cup\phi(A^*)\cup\phi(B^*)\\
&=\{s,t,u,v\}\cup\mathrm{cl}_{\beta Y}\bigg(\bigcup_{D\neq F\in\mathscr{C}(Y)}\{n_F\}\bigg)\\
&\cup\{s^*,t^*,u^*,v^*\}\cup\mathrm{cl}_{\beta Y}\bigg(\bigcup_{D^*\neq F\in\mathscr{C}(Y)}\{n_F\}\bigg)\\
&=\{s,t,u,v,n_D\}\cup\{s^*,t^*,u^*,v^*,n_{D^*}\}\cup\mathrm{cl}_{\beta Y}\bigg(\bigcup_{D,D^*\neq F\in\mathscr{C}(Y)}\{n_F\}\bigg).
\end{align*}
Thus, there are at least 9 minimal elements $Z$ in $\mathscr{Z}(Y)$ such that
\[Z\subseteq\phi(A\cup B\cup A^*\cup B^*),\]
namely,
\[Z=\{y\}\cup\{z\}\cup\mathrm{cl}_{\beta Y}\bigg(\bigcup_{D,D^*\neq F\in\mathscr{C}(Y)}\{n_F\}\bigg),\]
where $y$ is either $s$, $t$, $u$, $v$ or $n_D$ (from which at least three are distinct) and $z$ is either $s^*$, $t^*$, $u^*$, $v^*$ or $n_{D^*}$ (from which at least three are distinct). This contradiction shows that $D=D^*$.

Note that
\[A\cap B\cap A^*\cap B^*=\{c\}\cup\mathrm{cl}_{\beta X}\bigg(\bigcup_{C\neq E\in\mathscr{C}(X)}\{m_E\}\bigg),\]
and thus $A\cap B\cap A^*\cap B^*$ is a minimal element in $\mathscr{Z}(X)$. Therefore, its image
\[\phi(A)\cap\phi(B)\cap\phi(A^*)\cap\phi(B^*)\]
is a minimal element in $\mathscr{Z}(Y)$. Observe that
\[\phi(A)\cap\phi(B)\cap\mathrm{cl}_{\beta Y}D=\{s,t\}\cap\{u,v\}\]
and
\[\phi(A^*)\cap\phi(B^*)\cap\mathrm{cl}_{\beta Y}D=\{s^*,t^*\}\cap\{u^*,v^*\},\]
where $\{s,t\}\cap\{u,v\}$ and $\{s^*,t^*\}\cap\{u^*,v^*\}$ are both singletons. Since
\[\phi(A)\cap\phi(B)\cap\phi(A^*)\cap\phi(B^*)\cap\mathrm{cl}_{\beta Y}D\]
is non-empty, we have
\[\{s,t\}\cap\{u,v\}=\{s^*,t^*\}\cap\{u^*,v^*\}.\]
Therefore $\psi$ is a well defined mapping.

This completes the proof of the claim.
\end{proofoftheclaim}

\begin{xclaim}\label{HG}
$\psi$ is bijective.
\end{xclaim}

\begin{proofoftheclaim}
Note that $\phi^{-1}:\mathscr{Z}(Y)\rightarrow\mathscr{Z}(X)$ is an order-isomorphism. This induces a mapping
\[\psi^*:\bigcup_{D\in\mathscr{C}(Y)}(\mathrm{cl}_{\beta Y}D\setminus Y)\longrightarrow\bigcup_{C\in\mathscr{C}(X)}(\mathrm{cl}_{\beta X}C\setminus X)\]
which is defined analogous to $\psi$. Indeed, for an element $d$ of $\mathrm{cl}_{\beta Y}D\setminus Y$, where $D$ is a component of $Y$, we choose distinct elements $s$ and $t$ in $\mathrm{cl}_{\beta Y}D\setminus Y$ both distinct from $d$. We let $N$ and $M$ be (the previously defined) minimal elements of $\mathscr{Z}(Y)$ and $\mathscr{Z}(X)$, respectively, which are related by $\phi^{-1}(N)=M$. We let
\[S=\{s,d\}\cup\mathrm{cl}_{\beta Y}\bigg(\bigcup_{D\neq F\in\mathscr{C}(Y)}\{n_F\}\bigg)\quad\mbox{and}\quad T=\{t,d\}\cup\mathrm{cl}_{\beta Y}\bigg(\bigcup_{D\neq F\in\mathscr{C}(Y)}\{n_F\}\bigg),\]
and we define $\psi^*(d)$ to be the single element in the intersection $\{a,b\}\cap\{p,q\}$, where
\[\phi^{-1}(S)=\{a,b\}\cup\mathrm{cl}_{\beta X}\bigg(\bigcup_{C\neq E\in\mathscr{C}(X)}\{m_E\}\bigg)\]
and
\[\phi^{-1}(T)=\{p,q\}\cup\mathrm{cl}_{\beta X}\bigg(\bigcup_{C\neq E\in\mathscr{C}(X)}\{m_E\}\bigg),\]
$C$ is a component of $X$, and $a$, $b$, $p$ and $q$ are elements of $\mathrm{cl}_{\beta X}C\setminus X$. Observe that, if we let $A=\phi^{-1}(S)$ and $B=\phi^{-1}(T)$, then $S=\phi(A)$ and $T=\phi(B)$, and therefore, by the definition of $\psi$ we have
\[\psi\big(\psi^*(d)\big)=d.\]
In particular, this shows that $\psi$ is surjective. A similar argument shows that
\[\psi^*\big(\psi(c)\big)=c\]
for every element $c$ in $\mathrm{cl}_{\beta X}C\setminus X$, where $C$ is a component of $X$. Therefore, $\psi$ is also injective. In particular,
\[\psi^{-1}=\psi^*.\]

This completes the proof of the claim.
\end{proofoftheclaim}

\begin{xclaim}\label{PHFSD}
For every component $C$ of $X$ there is a unique component $C_\psi$ of $Y$ such that
\[\psi(\mathrm{cl}_{\beta X}C\setminus X)=\mathrm{cl}_{\beta Y}C_\psi\setminus Y.\]
\end{xclaim}

\begin{proofoftheclaim}
Let $c$ be an element in $\mathrm{cl}_{\beta X}C\setminus X$ where $C$ is a component of $X$. Let $a$ be an element in $\mathrm{cl}_{\beta X}C\setminus X$ distinct from $c$. Choose some $b$ in $\mathrm{cl}_{\beta X}C\setminus X$ distinct from both $a$ and $c$. Let
\[S=\{a,b\}\cup\mathrm{cl}_{\beta X}\bigg(\bigcup_{C\neq E\in\mathscr{C}(X)}\{m_E\}\bigg),\quad T=\{b,c\}\cup\mathrm{cl}_{\beta X}\bigg(\bigcup_{C\neq E\in\mathscr{C}(X)}\{m_E\}\bigg)\]
and
\[U=\{c,a\}\cup\mathrm{cl}_{\beta X}\bigg(\bigcup_{C\neq E\in\mathscr{C}(X)}\{m_E\}\bigg).\]
Then $\psi(a)$ and $\psi(c)$, by definition, are both in $\mathrm{cl}_{\beta Y}D\setminus Y$ where $D$ is the (unique) component of $Y$ such that
\[\big|\phi(U)\cap(\mathrm{cl}_{\beta Y}D\setminus Y)\big|>1.\]
Define
\[C_\psi=D.\]
Therefore
\begin{equation}\label{OPJIS}
\psi(\mathrm{cl}_{\beta X}C\setminus X)\subseteq\mathrm{cl}_{\beta Y}C_\psi\setminus Y.
\end{equation}
Symmetry (using the order-isomorphism $\phi^{-1}$ instead of $\phi$) implies that
\[\psi^{-1}(\mathrm{cl}_{\beta Y}C_\psi\setminus Y)\subseteq\mathrm{cl}_{\beta X}(C_\psi)_{\psi^{-1}}\setminus X.\]
Therefore
\begin{equation}\label{OD}
\mathrm{cl}_{\beta Y}C_\psi\setminus Y=\psi\big(\psi^{-1}(\mathrm{cl}_{\beta Y}C_\psi\setminus Y)\big)\subseteq\psi\big(\mathrm{cl}_{\beta X}(C_\psi)_{\psi^{-1}}\setminus X\big)
\end{equation}
which combined with (\ref{OPJIS}) implies that
\[\psi(\mathrm{cl}_{\beta X}C\setminus X)\subseteq\psi\big(\mathrm{cl}_{\beta X}(C_\psi)_{\psi^{-1}}\setminus X\big).\]
Therefore
\[\mathrm{cl}_{\beta X}C\setminus X\subseteq\mathrm{cl}_{\beta X}(C_\psi)_{\psi^{-1}}\setminus X.\]
But, the closure in $\beta X$ of components of $X$ are pairwise disjoint (as components of $X$ are pairwise disjoint open and closed subspaces of $X$), thus
\begin{equation}\label{PDF}
(C_\psi)_{\psi^{-1}}=C.
\end{equation}
Combining (\ref{OD}) and (\ref{PDF}) proves the reverse inclusion in (\ref{OPJIS}).

This completes the proof of the claim.
\end{proofoftheclaim}

Observe that in the proof of the above claim we have indeed defined a mapping from the set of all components of $X$ to the set of all components of $Y$ assigning the component $C$ of $X$ to the component $C_\psi$ of $Y$. By (\ref{PDF}) this mapping is injective. That this mapping is also surjective follows from a similar argument which shows that
\[(D_{\psi^{-1}})_\psi=D\]
for every component $D$ of $Y$.

\begin{xclaim}
For every component $C$ of $X$ the restricted mapping
\[\psi:\mathrm{cl}_{\beta X}C\setminus X\longrightarrow\mathrm{cl}_{\beta Y}C_\psi\setminus Y\]
is a homeomorphism.
\end{xclaim}

\begin{proofoftheclaim}
Let $C$ be a component of $X$. By the previous claim we need to show only that the restricted mapping $\psi$ and its inverse are continuous. Note that from the definition of $\psi$ we have
\[\phi\bigg(\{c\}\cup\mathrm{cl}_{\beta X}\bigg(\bigcup_{C\neq E\in\mathscr{C}(X)}\{m_E\}\bigg)\bigg)=\big\{\psi(c)\big\}\cup\mathrm{cl}_{\beta Y}\bigg(\bigcup_{C_\psi\neq F\in\mathscr{C}(Y)}\{n_F\}\bigg)\]
for every $c$ in $\mathrm{cl}_{\beta X}C\setminus X$.

Let $K$ be a closed subspace of $\mathrm{cl}_{\beta X}C\setminus X$. Let
\[Z_k=\{k\}\cup\mathrm{cl}_{\beta X}\bigg(\bigcup_{C\neq E\in\mathscr{C}(X)}\{m_E\}\bigg)\]
for every $k$ in $K$. By (the proof of) Theorem \ref{KDSG} we have
\[\phi\bigg(\mathrm{cl}_{\beta X}\bigg(\bigcup_{k\in K}Z_k\bigg)\bigg)=\phi\bigg(\bigvee_{k\in K}Z_k\bigg)=\bigvee_{k\in K}\phi(Z_k)=\mathrm{cl}_{\beta Y}\bigg(\bigcup_{k\in K}\phi(Z_k)\bigg),\]
where the first and the second $\vee$ are to denote the least upper bounds in $\mathscr{Z}(X)$ and $\mathscr{Z}(Y)$, respectively. But
\[\bigcup_{k\in K}Z_k=K\cup\mathrm{cl}_{\beta X}\bigg(\bigcup_{C\neq E\in\mathscr{C}(X)}\{m_E\}\bigg)\]
and
\[\bigcup_{k\in K}\phi(Z_k)=\psi(K)\cup\mathrm{cl}_{\beta Y}\bigg(\bigcup_{C_\psi\neq F\in\mathscr{C}(Y)}\{n_F\}\bigg).\]
Therefore, from the above relations we have
\begin{align}
\phi\bigg(K\cup\mathrm{cl}_{\beta X}\bigg(\bigcup_{C\neq E\in\mathscr{C}(X)}\{m_E\}\bigg)\bigg)=\mathrm{cl}_{\beta Y}\psi(K)\cup\mathrm{cl}_{\beta Y}\bigg(\bigcup_{C_\psi\neq F\in\mathscr{C}(Y)}\{n_F\}\bigg).
\label{RSDS}
\end{align}
Note that $\mathrm{cl}_{\beta Y}C_\psi\setminus Y$ contains $\psi(K)$ by the previous claim, and therefore contains its closure in $\beta Y$. Let $L=\mathrm{cl}_{\beta Y}\psi(K)$. A similar argument shows that
\[\phi^{-1}\bigg(L\cup\mathrm{cl}_{\beta Y}\bigg(\bigcup_{C_\psi\neq F\in\mathscr{C}(Y)}\{n_F\}\bigg)\bigg)=\mathrm{cl}_{\beta X}\psi^{-1}(L)\cup\mathrm{cl}_{\beta X}\bigg(\bigcup_{(C_\psi)_{\psi^{-1}}\neq E\in\mathscr{C}(X)}\{m_E\}\bigg)\]
with
\[\psi^{-1}(L)\subseteq\mathrm{cl}_{\beta X}(C_\psi)_{\psi^{-1}}\setminus X.\]
Observe that $(C_\psi)_{\psi^{-1}}=C$. Therefore
\begin{align}
\phi\bigg(\mathrm{cl}_{\beta X}\psi^{-1}(L)\cup\mathrm{cl}_{\beta X}\bigg(\bigcup_{C\neq E\in\mathscr{C}(X)}\{m_E\}\bigg)\bigg)=L\cup\mathrm{cl}_{\beta Y}\bigg(\bigcup_{C_\psi\neq F\in\mathscr{C}(Y)}\{n_F\}\bigg).
\label{FSD}
\end{align}
Comparing (\ref{RSDS}) and (\ref{FSD}) it follows that $K=\mathrm{cl}_{\beta X}\psi^{-1}(L)$. In particular, then $\psi^{-1}(L)\subseteq K$, or, equivalently, $L\subseteq\psi(K)$. On the other hand, $\psi(K)\subseteq L$ by the definition of $L$. Therefore $\psi(K)=L$ and thus $\psi(K)$ is compact and therefore closed in $\mathrm{cl}_{\beta Y}C_\psi\setminus Y$.

A similar argument shows that every closed subspace of $\mathrm{cl}_{\beta Y}C_\psi\setminus Y$ is mapped by $\psi^{-1}$ onto a closed subspace of $\mathrm{cl}_{\beta X}C\setminus X$.

This completes the proof of the claim.
\end{proofoftheclaim}

The proof of the theorem is now complete.
\end{proof*}

\subsubsection*{Acknowledgements}

The author wishes to thank Professor Micha{\l} Kukie{\l}a for his comment on an early draft of this article.

\end{document}